\documentclass{amsart}

\usepackage{graphicx}
\usepackage{amsfonts}
\usepackage{amssymb}
\usepackage{amsthm}
\usepackage{amsmath}
\usepackage[all]{xy}

\theoremstyle{plain}
\newtheorem{theorem}[subsection]{Theorem}
\newtheorem{lemma}[subsection]{Lemma}
\newtheorem{proposition}[subsection]{Proposition}
\newtheorem{corollary}[subsection]{Corollary}

\theoremstyle{definition}

\newtheorem{definition}[subsection]{Definition}

\newtheorem{remark}[subsection]{Remark}
\newtheorem{notation}[subsection]{Notation}
\newtheorem{example}[subsection]{Example}

\newcommand{\comp}{\circ}
\newcommand{\defn}{\textbf}

\renewcommand{\to}{\longrightarrow}
\newcommand{\To}{\Longrightarrow}
\newcommand{\del}{\partial}
\newcommand{\eps}{\epsilon}
\newcommand{\iso}{\cong}
\newcommand{\im}{\ensuremath{\mathsf{im\,}}}
\renewcommand{\Im}{\ensuremath{\mathrm{Im}}}
\newcommand{\coker}{\ensuremath{\mathsf{coker\,}}}
\newcommand{\Cok}{\ensuremath{\mathrm{Q}}}
\renewcommand{\ker}{\ensuremath{\mathsf{ker\,}}}
\newcommand{\Ker}{\ensuremath{\mathrm{K}}}

\renewcommand{\H}{\ensuremath{\mathrm{H}}}
\newcommand{\Zu}{\ensuremath{\mathrm{Z}}}
\newcommand{\Nu}{\ensuremath{\mathrm{N}}}
\newcommand{\Ac}{\ensuremath{\mathcal{A}}}

\newcommand{\Cc}{\ensuremath{\mathcal{C}}}
\newcommand{\Dc}{\ensuremath{\mathcal{D}}}
\newcommand{\Ec}{\ensuremath{\mathcal{E}}}
\newcommand{\Pc}{\ensuremath{\mathcal{P}}}
\newcommand{\Sc}{\ensuremath{\mathcal{S}}}

\newcommand{\Rng}{\ensuremath{\mathsf{Rng}}}
\newcommand{\Comm}{\ensuremath{\mathsf{Comm}}}

\newcommand{\XMod}{\ensuremath{\mathsf{XMod}}}

\newcommand{\Alg}{\ensuremath{\mathsf{Alg}}}
\newcommand{\G}{\ensuremath{\mathbb{G}}}
\newcommand{\K}{\ensuremath{\mathbb{K}}}

\newcommand{\Z}{\ensuremath{\mathbb{Z}}}
\newcommand{\N}{\ensuremath{\mathbb{N}}}

\newcommand{\Set}{\ensuremath{\mathsf{Set}}}

\newcommand{\Gp}{\ensuremath{\mathsf{Gp}}}

\newcommand{\Mod}{\ensuremath{\mathsf{Mod}}}

\newcommand{\Hom}{\ensuremath{\mathrm{Hom}}}
\newcommand{\Tor}{\ensuremath{\mathrm{Tor}}}

\renewcommand{\hom}{\ensuremath{\mathrm{Hom}}}

\newcommand{\pr}{\ensuremath{\mathrm{pr}}}
\newcommand{\noproof}{\hfill \qed}
\newcommand{\Pepi}{\ensuremath\text{${\mathcal{P}}$-epi}}
\newcommand{\Eproj}{\ensuremath\text{${\mathcal{E}}$-proj}}

\newbox\skewpullbackbox
\setbox\skewpullbackbox=\hbox{\xy 0;<1mm,0mm>: \POS(4,0)\ar@{-} (-4,0) \ar@{-} (8,4)
\endxy}

\newbox\ksewpullbackbox
\setbox\ksewpullbackbox=\hbox{\xy 0;<1mm,0mm>: \POS(0,-8)\ar@{-} (0,-4) \ar@{-} (4,-4)
\endxy}

\newbox\pullbackbox
\setbox\pullbackbox=\hbox{\xy 0;<1mm,0mm>: \POS(4,0)\ar@{-} (0,0) \ar@{-} (4,4)
\endxy}

\newbox\pushoutbox
\setbox\pushoutbox=\hbox{\xy 0;<1mm,0mm>: \POS(0,4)\ar@{-} (0,0) \ar@{-} (4,4)
\endxy}

\begin{document}

\newdir{>>}{{}*!/3.5pt/:(1,-.2)@^{>}*!/3.5pt/:(1,+.2)@_{>}*!/7pt/:(1,-.2)@^{>}*!/7pt/:(1,+.2)@_{>}}
\newdir{ >>}{{}*!/8pt/@{|}*!/3.5pt/:(1,-.2)@^{>}*!/3.5pt/:(1,+.2)@_{>}}
\newdir{ |>}{{}*!/-3.5pt/@{|}*!/-8pt/:(1,-.2)@^{>}*!/-8pt/:(1,+.2)@_{>}}
\newdir{ >}{{}*!/-8pt/@{>}}
\newdir{>}{{}*:(1,-.2)@^{>}*:(1,+.2)@_{>}}
\newdir{<}{{}*:(1,+.2)@^{<}*:(1,-.2)@_{<}}

\title{A Comparison Theorem for Simplicial Resolutions}

\author{Julia Goedecke}
\email{julia.goedecke@cantab.net}
\address{DPMMS\\
University of Cambridge\\
Wilberforce Road\\
Cambridge CB3 0WB\\
United Kingdom}

\author{Tim Van der Linden}
\email{tvdlinde@vub.ac.be}
\address{Vakgroep Wiskunde\\
Vrije Universiteit Brussel\\
Pleinlaan~2\\
1050~Brussel\\
Belgium}

\thanks{The second author's research was supported by York University, Toronto. Published as: J.~Goedecke and T.~Van~der Linden, \emph{A comparison theorem for simplicial resolutions}, J.~Homotopy and Related Structures \textbf{2} (2007), no.~1, 109--126.}



\begin{abstract}
It is well known that Barr and Beck's definition of comonadic homology makes sense also with a functor of coefficients taking values in a semi-abelian category instead of an abelian one. The question arises whether such a homology theory has the same convenient properties as in the abelian case. Here we focus on independence of the chosen comonad: conditions for homology to depend on the induced class of projectives only.
\end{abstract}




\maketitle

\section*{Introduction}

Given two comonads $\G$ and $\K$ on an arbitrary category $\Cc$, when do they induce the same comonadic homology theory? In their paper~\cite{Barr-Beck}, Barr and Beck gave a sufficient condition for this to happen: when $\G$ and $\K$ generate the same class of projective objects. More precisely, they showed that $\H_n(-,E)_{\G}\iso \H_n(-,E)_{\K}$ as functors $\Cc\to\Ac$, for any $n\in \N$ and for any functor $E\colon{\Cc\to\Ac}$ to an abelian category $\Ac$. 

With this paper, we extend their result to the semi-abelian case, i.e.\ the situation where the homology theory takes coefficients in a functor $E\colon{\Cc\to\Ac}$ to a semi-abelian category $\Ac$. To obtain the same conclusion, we need a condition on $\Cc$ (or on the comonads $\G$ and $\K$): homming from a ($\G$- or $\K$-) projective object into a ($\G$- or $\K$-) simplicial resolution gives a Kan simplicial set. 

The examples reveal that this condition is not unreasonably strong. For instance, if $\Cc$ is additive then homming from \emph{any} object into \emph{any} simplicial object gives a Kan simplicial set, and if $\Cc$ is a regular Mal'tsev category---for instance, $\Cc$ could be semi-abelian---then the requirement holds as soon as the $\G$-projective objects are also regular projectives. Further examples are given in Section~\ref{Section-Examples}.

The technique we use to obtain our result is necessarily different from the one employed by Barr and Beck. Since their coefficient functors take values in an abelian category, they can extend the homology theory on $\Cc$ to a homology theory on the free additive category over $\Cc$; then this additive structure takes part in the proof that $\H_n(-,E)_{\G}$ depends only on the class of $\G$-projectives. As a semi-abelian category is only additive when it is abelian, this approach is bound to fail in the semi-abelian case.

However, our strategy is still simple and straightforward: we show that, given a projective class $\Pc$, any two $\Pc$-resolutions of an object $X$ are homotopically equivalent---hence they have the same homology. The advantage of our method is that it shows that \emph{any} $\Pc$-resolution will give the same homology of $X$; thus it is possible to use resolutions not coming from a comonad, if this turns out to be easier. The only subtlety is in the definition of a \defn{$\Pc$-resolution} of $X$: this is an augmented simplicial object $A=(A_n)_{n\geq-1}$ such that $A_{-1}=X$, $A_n\in \Pc$ for $n\geq0$, and for every object $P\in \Pc$ the augmented simplicial set $\Hom(P,A)$ is Kan and contractible. It is here that the relative Kan condition on $\Cc$ emerges: we have to know whether, when $\Pc$ is the class of $\G$-projectives, a $\G$-resolution is always a $\Pc$-resolution. 

In Section~\ref{Section-Semi-abelian-Context} we give an overview of the required results from the theory of semi-abelian categories, and the basics on homology in this context. The definition of a $\Pc$-resolution is further explained in Section~\ref{Section-Resolutions}. Section~\ref{Section-Comparison-Theorem} is devoted to the Comparison Theorem~\ref{comparison}: if $P$ is a simplicial object over $X$ with each $P_i\in\Pc$, and $A$ is a simplicial object over $Y$ where all augmented simplicial sets $\Hom(P_i,A)$ are contractible and Kan, then any map $f\colon {X\to Y}$ extends to a semi-simplicial map $f\colon {P\to A}$, and any two such extensions are simplicially homotopic. In this section we also relate our comparison theorem to that of Tierney and Vogel \cite{Tierney-Vogel2}, which uses a different definition of resolution, in a category with finite limits.

The Comparison Theorem is used in Section~\ref{Section-Main-Theorem} to prove the main point of our paper, Theorem~\ref{Main-Theorem}: under the condition on $\Cc$ mentioned above, any two comonads $\G$ and $\K$ that generate the same class of projectives induce isomorphic homology theories. We obtain it as an immediate consequence of Corollary~\ref{Corollary-Homotopic-then-Same-Homology} which states that, in a semi-abelian category, simplicially homotopic maps have the same homology: if $f\simeq g$ then, for any $n\in \N$, $\H_{n}f= \H_{n}g$.

\section{The semi-abelian context}\label{Section-Semi-abelian-Context}

\subsection{Semi-abelian categories}
It is well known that the homological diagram lemmas like the Five Lemma and the Snake Lemma can be proved in any abelian category, but are also true in the category of groups, and in other categories `sufficiently close' to it. Semi-abelian categories were invented to capture this more general context. 

\begin{definition} A category $\Ac$ is called \defn{Barr exact} \cite{Barr} when it is \defn{regular} (i.e.\ finitely complete with coequalizers of kernel pairs and pullback-stable regular epimorphisms) and every equivalence relation in $\Ac$ is a kernel pair. 

A pointed and regular category $\Ac$ is called \defn{Bourn protomodular} \cite{Bourn1991} when the \defn{(regular) Short Five Lemma} holds: 
for every commutative diagram
$$\xymatrix{K[f'] \ar@{{ |>}->}[r]^-{\ker f'} \ar[d]_-u & X' \ar@{-{ >>}}[r]^-{f'} \ar[d]^-v & Y' \ar[d]^-w \\ 
K[f] \ar@{{ |>}->}[r]_-{\ker f} & X \ar@{-{ >>}}[r]_-{f} & Y}$$
such that $f$ and $f'$ are regular epimorphisms, $u$ and $w$ being isomorphisms implies that $v$ is an isomorphism.
\end{definition}

\begin{definition}\label{def-semi-abelian}\cite{Janelidze-Marki-Tholen}
A \defn{semi-abelian category} is a category $\Ac$ which is pointed, Barr exact and Bourn protomodular, and has binary coproducts.
\end{definition}

All finite limits and colimits exist in a semi-abelian category, and as $\Ac$ is pointed, there is a zero object: an initial object which is also terminal. Thus we can form kernels and cokernels. The \defn{kernel} $\ker{f}\colon {\Ker[f]\to X}$ of a map $f\colon {X\to Y}$ is the pullback of $0\to Y$ along $f$, and dually the \defn{cokernel} $\coker{f}\colon {Y\to \Cok[f]}$ is the pushout of $X\to 0$ along $f$. It can be shown that $\Ker[f]=0$ iff $f$ is a monomorphism, and $\Cok [f]=0$ iff $f$ is a regular epimorphism. Also, any regular epimorphism is the cokernel of its kernel.

The zero object also gives us \defn{zero maps}, which are maps that factor through $0$.

\begin{example} 
Any abelian category is semi-abelian, as is the category $\Gp$ of groups. More generally, any variety of $\Omega$-groups is semi-abelian: for instance, the categories of non-unital rings, Lie algebras, (pre)crossed modules. 
\end{example}

As hinted above, the basic homological diagram lemmas can be proved in any semi-abelian category (see \cite{Borceux-Semiab,Borceux-Bourn}). But many other things that are true in an abelian category do not hold in a semi-abelian one: the $\Hom$-sets need not be abelian groups, binary products do not coincide with binary coproducts in general, and maps cannot in general be factored into a cokernel followed by a kernel. But as a semi-abelian category is \emph{regular}, we can factor any map $f\colon {X\to Y}$ into a regular epimorphism (in fact, a cokernel) $X\to \Im[f]$ followed by a monomorphism $\im{f}\colon {\Im[f]\to Y}$. This monomorphism is called the \defn{image} of $f$. If it is a normal monomorphism, i.e.\ the kernel of another map, $f$ is called \defn{proper}. This factorisation is unique up to isomorphism, and it allows us to define exact sequences in a semi-abelian category.

\begin{definition}\label{def-exact-sequences}
A sequence of morphisms 
$$
\xymatrix{X \ar[r]^-{f} & Y \ar[r]^-g & Z}
$$
is called \defn{exact (at $Y$)} if $\im{f}=\ker{g}$.
\end{definition}

A sequence $0\to X\to Y$ is exact iff the map $X\to Y$ is a monomorphism, and $X\to Y\to 0$ is exact iff $X\to Y$ is a regular epimorphism.  Note that a non-proper map can never occur as the first map of an exact sequence. 

\subsection{Homology in semi-abelian categories}
As usual a \defn{chain complex} is a sequence of maps $(d_{n}\colon {C_n \to C_{n-1}})_{n\in\Z}$ with $d_{n-1}\comp d_n=0$ for all $n$. A chain complex is called \defn{proper} when all the maps $d_n$ are proper maps. 

\begin{definition}\label{def-homology}\cite{EverVdL2}
Let $C$ be a proper chain complex in a semi-abelian category. The $n$-th homology object $\H_n C$ is the cokernel of $d'_{n+1}\colon {C_{n+1}\to \Ker[d_{n}]}$, the factorisation of $d_{n+1}\colon{C_{n+1}\to C_{n}}$ over $\ker d_{n}$.
\end{definition}

It is easy to see from the definition that $C$ is exact at $C_n$ if and only if $\H_nC=0$. Thus the homology gives us the usual detection of exactness. When talking about homology, we only consider proper chain complexes, because otherwise this property is false. Consider, for instance, the following example in the category of groups. In $\Gp$ a monomorphism is normal iff it is the inclusion of a normal subgroup. Define a chain complex by taking $d_1$ to be the inclusion of $A_4$ into $A_5$, and all other objects to be zero. Since $A_5$ is simple, $d_1$ is not proper, and all objects $\H_nC$ are zero, but clearly $C$ is not exact at $C_0$.\\

Given a simplicial object $A=(A_n)_{n\geq 0}$ (with face operators $\del_{i}\colon {A_{n}\to A_{n-1}}$ for $i\in [n]=\{0,\dots ,n \}$ and $n\in \N_{>0}$, and degeneracy operators $\sigma_{i}\colon {A_{n}\to A_{n+1}}$, for $i\in [n]$ and $n\in \N$, subject to the simplicial identities) in a semi-abelian category $\Ac$, we can define the homology of $A$ by going via the \emph{Moore complex} of $A$. 

\begin{definition}\label{def-Moore-complex}
Let $A$ be a simplicial object in a semi-abelian category $\Ac$. The \defn{Moore complex} $\Nu(A)$ has as objects $\Nu_0A=A_0$, $\Nu_{-n}A=0$ and 
$$\Nu_nA=\bigcap^{n-1}_{i=0}\Ker[\del_i\colon {A_n\to A_{n-1}}] = \Ker[(\del_i)_{i\in [n-1]}\colon {A_n\to A^n_{n-1}}],$$
for $n\geq1$, and boundary maps 
$d_n=\del_n\comp \bigcap_{i}\ker{\del_i}\colon {\Nu_nA\to \Nu_{n-1}A}$
 for $n\geq1$.
 
The \defn{object of $n$-cycles} is $\Zu_nA=\Ker[d_n]=\bigcap^{n}_{i=0}\Ker[\del_i\colon {A_n\to A_{n-1}}]$ for $n\geq1$. We write $\Zu_0A=A_0$.
\end{definition}

The Moore complex of a simplicial object is always a proper chain complex \cite[Theorem~3.6]{EverVdL2}; thus we can define 
$$\H_n A=\H_n \Nu(A)$$
for a simplicial object $A$. In the abelian case, the homology of the Moore complex is the same as the homology of the \defn{unnormalised} chain complex $\mathrm{C}(A)$ of $A$, where $\mathrm{C}_nA=A_n$ and $d_n=\del_0-\del_1+\cdots+(-1)^n\del_n$.  

Notice that the Moore complex and thus the homology of a simplicial object only involve the face maps $\del_i$, and not the degeneracies $\sigma_i$. Thus in this context it is enough to consider semi-simplicial maps between simplicial objects, i.e.\ maps that commute with the $\del_i$ but not necessarily with the $\sigma_i$. 

\subsection{Comonadic homology}
Let $\Cc$ be an arbitrary category and
\[
\G =(G\colon \Cc\to \Cc,\quad \delta\colon G\To G^{2} ,\quad \epsilon\colon G\To 1_{\Cc})
\]
a comonad on $\Cc$. Given an object $X$ in $\Cc$, $\G$ gives rise to an augmented simplicial object over $X$, where the objects are $G^{n+1}X$, and the maps are $\del_i=G^i\eps_{G^{n-i}X}\colon {G^{n+1}X\to G^nX}$ and $\sigma_i=G^i\delta_{G^{n-i}X}\colon {G^{n+1}X\to G^{n+2}X}$. We write $\G X$ for this simplicial object.
\[
\xymatrix{{\cdots} \ar@<1.5ex>[r] \ar@<.5ex>[r] \ar@<-.5ex>[r] \ar@<-1.5ex>[r] & G^{3}X \ar[r] \ar@<1ex>[r] \ar@<-1ex>[r] & G^{2}X \ar@<.5ex>[r] \ar@<-.5ex>[r] & GX \ar[r] & X}
\]
Given a functor $E\colon {\Cc\to \Ac}$ into a semi-abelian category, it is then possible to form the homology of the resulting simplicial object $E\G X$. This naturally generalises Barr-Beck comonadic homology~\cite{Barr-Beck} to the semi-abelian context:

\begin{definition}\cite{EverVdL2}\label{Definition-Homology}
Let $\Cc$ be a category equipped with a comonad $\G$ and $\Ac$ a semi-abelian category. Let $E\colon \Cc\to \Ac$ be a functor. For $n\geq 1$, the object
\[
\H_{n} (X,E)_{\G}=\H_{n-1} \Nu E \G X
\]
is the \defn{$n$-th homology object of $X$ (with coefficients in $E$) relative to the comonad~$\G$}. This defines a functor $\H_{n}(-,E)_\G \colon \Cc\to \Ac$, for every $n\geq 1$. 
\end{definition}

The dimension shift here is not present in Barr and Beck's original definition, but was introduced in~\cite{EverVdL2} to make it better adjusted to the non-abelian examples (homology of groups, Lie algebras, crossed modules) which traditionally have a shifted numbering. Definition~\ref{Definition-Homology} is consistent with the Hopf formulae~\cite{EGVdL} which exist in the non-abelian case.

\subsection{Technical lemmas}
We now give some technical lemmas that hold in our context, which we will need later on. The next lemma is true already in regular categories, and is taken from~\cite{Carboni-Kelly-Pedicchio}. It will be used to show that a particular morphism is regular epic.

\begin{lemma}\label{Lemma-show-regular-epi}
Let $\Ac$ be a regular category. A map $y\colon {Z_{0}\to Y}$ factorises through the image of a map $f\colon {X\to Y}$ if and only if there is a regular epimorphism $z\colon {Z\to Z_{0}}$ and a map $x\colon {Z\to X}$ with $yz=fx$.\noproof
\end{lemma}

\begin{remark}\label{Remark-show-regular-epi}
Of course, if we can show that \emph{every} map $y\colon {Z_0\to Y}$ factorises through the image of a given map $f\colon {X\to Y}$, this map $f$ is a regular epimorphism. Hence we can use Lemma~\ref{Lemma-show-regular-epi} to show that a map $f\colon {X\to Y}$ is regular epic by proving that, for every map $y\colon {Z_{0}\to Y}$, there is a regular epimorphism $z\colon {Z\to Z_{0}}$ and a map $x\colon {Z\to X}$ with $yz=fx$. Thus an argument which would otherwise involve the use of projective objects (proving that $f$ is a regular epimorphism by showing that every map $y$ with a projective domain $Z_{0}$ factors over $f$) can be replaced by an argument which uses regular epimorphisms only, and thus the requirement that enough projectives exist may be avoided. We use this in the proof of Proposition~\ref{Proposition-phi_0-Homology-Isomorphism}.
\end{remark}

\begin{lemma}\label{Lemma-pushout-implies-cokernels-iso} Let $\Ac$ be a pointed category. Consider the diagram 
$$
\xymatrix{X' \ar[r]^-{f'} \ar[d] & Y' \ar@{-{ >>}}[r] \ar[d] & \Cok[f'] \ar[d]\\
X \ar[r]_-{f} & Y \ar@{-{ >>}}[r] & \Cok[f]}
$$
where $\Cok[f]$ and $\Cok[f']$ are the cokernels of $f$ and $f'$ respectively. If the left square is a pushout, then the induced map between the cokernels is an isomorphism.\noproof
\end{lemma}

The next definition is slightly more general than the usual definition of a regular pushout, where all the maps are demanded to be regular epimorphims \cite{Bourn2003, Carboni-Kelly-Pedicchio}. We will need it in this more general form.

\begin{definition}\label{Definition-regepi-Pushout}
Let $\Ac$ be a semi-abelian category. A square in $\Ac$ with horizontal regular epimorphisms
$$
\xymatrix{ X' \ar@{-{ >>}}[r]^-{f'} \ar[d]_-{x} & Y' \ar[d]^-y \\
X \ar@{-{ >>}}[r]_-f & Y}
$$
is called a \defn{regular pushout} when the comparison map $(x,f')\colon  {X'\to X\times_{Y}Y'}$ to the pullback $X\times_{Y}Y'$ of $y$ along $f$ is a regular epimorphism. (The maps $x$ and $y$ are \emph{not} demanded to be regular epimorphisms.)
\end{definition}
A regular pushout is always a pushout, but a pushout need not be a regular pushout.

\section{Simplicial resolutions}\label{Section-Resolutions}

\subsection{Projective classes} We need to consider simplicial resolutions relative to a chosen class of projectives. Here we recall the definition of a \emph{projective class} and give some examples.

\begin{definition}\label{Definition-Projective-Class}
Let $\Cc$ be a category, $P$ an object and $e\colon {X\to Y}$ a morphism of $\Cc$. Then $P$ is called \defn{$e$-projective}, and $e$ is called \defn{$P$-epic}, if the induced map
$$\hom (P,e)=e\comp (\cdot)\colon {\hom (P,X)\to \hom (P,Y)}$$
is a surjection. Let $\Pc$ be a class of objects of $\Cc$. A morphism $e$ is called \defn{$\Pc$-epic} if it is $P$-epic for every $P\in \Pc$; the class of all $\Pc$-epimorphisms is denoted by $\Pepi$. Let $\Ec$ be a class of morphisms in $\Cc$. An object $P$ is called \defn{$\Ec$-projective} if it is $e$-projective for every $e$ in $\Ec$; the class of $\Ec$-projective objects is denoted $\Eproj$. $\Cc$ is said \defn{to have enough $\Ec$-projectives} if for every object $Y$ there is a morphism $P\to Y$ in $\Ec$ with $P$ in $\Eproj$.

A \defn{projective class on $\Cc$} is a pair $(\Pc ,\Ec)$, $\Pc$ a class of objects of $\Cc$, $\Ec$ a class of morphisms of $\Cc$, such that $\Pc =\Eproj$, $\Pepi=\Ec$ and $\Cc$ has enough $\Ec$-projectives. Since, given a projective class $(\Pc ,\Ec)$, $\Pc$ and $\Ec$ determine each other, we will sometimes abusively write \defn{the projective class $\Pc$} or \defn{the projective class $\Ec$}.
\end{definition} 

It is easy to see that any retract of a projective object is also projective, as is any coproduct of projectives. 

\begin{example}
If $\Ec$ is the class of regular epimorphisms, $\Pc$ is called the class of \defn{regular projectives}. In a variety, the class of regular projectives is generated by the free objects, hence there are enough projectives. 
\end{example}

\subsection{The projective class generated by a comonad} The regular projectives in a variety $\Cc$ are also generated by the values of the canonical comonad $\Cc$, induced by the forgetful functor to $\Set$. More generally, \emph{any} comonad on a category $\Cc$ generates a projective class:

\begin{definition} Let $\G=(G,\eps,\delta)$ be a comonad on a category $\Cc$.  An object $P$ in $\Cc$ is called \defn{$\G$-projective} if it is in the projective class $(\Pc_{\G},\Ec_{\G})$ generated by the objects of the form $GY$. A map in $\Ec_{\G}$ is called a \defn{$\G$-epimorphism}.
\end{definition}

The $X$-component $\eps_X\colon {GX\to X}$ of the counit $\eps$ is always a $\G$-epimorphism. Indeed, any map $f\colon {GY\to X}$ factors over $\eps_{X}$ as $Gf\comp \delta_{Y}$, because  $\eps_{X}\comp Gf\comp \delta_{Y} = f\comp \eps_{GY}\comp \delta_{Y}=f$. It is now clear that $\Cc$ has enough projectives of this class, since for any $X$ we have $\eps_X\colon {GX\to X}$. 

This definition coincides with the definition of $\G$-projectives in~\cite{Barr-Beck}. There a \emph{$\G$-projective object} is an object $X$ which admits a map $s\colon {X\to GX}$ such that $\eps_X s=1_X$. Indeed, if $X\in \Pc$, then the identity on $X$ factors over the $\Pc$-epimorphism $\eps_X$, which gives the splitting $s$.

\subsection{The relative Kan property}\label{Section-Relative-Kan-Property}

A classical technical property simplicial sets may have is the Kan property. Kan simplicial sets are exactly the fibrant ones (in the usual model structure on $\Sc\Set$) and may be described as follows. Let $S$ be a simplicial set and $n\geq 1$, $k\in [n]$ natural numbers. An \defn{$(n,k)$-horn} in $S$ is a sequence $(s_{i})_{i\in [n]\setminus k}$ of elements of $S_{n-1}$ satisfying $\del_{i}(s_{j})=\del_{j-1}(s_{i})$ for all $i<j$ and $i,j\neq k$. A \defn{filler} of an $(n,k)$-horn $(s_{i})_{i}$ is an element $s$ of $S_{n}$ satisfying $\del_{i}(s)=s_{i}$ for all $i\neq k$. A simplicial set $S$ is \defn{Kan} when every horn in $S$ has a filler.

We need the simplicial objects in the category $\Cc$ to satisfy a similar property, but relative to a chosen projective class $\Pc$ on $\Cc$. We will say that a simplicial object $A$ is \defn{Kan (relative to $\Pc$)} when for every object $P\in \Pc$ the simplicial set $\Hom(P,A)$ is Kan.

\begin{example}
If $\Cc$ is regular with enough regular projectives and $\Pc$ the induced projective class, saying that $A$ is Kan relative to $\Pc$ is the same as saying that the simplicial object $A$ is Kan, in the internal sense of~\cite{Carboni-Kelly-Pedicchio}. Every simplicial object of $\Cc$ has this Kan property if and only if $\Cc$ is a Mal'tsev category~\cite[Theorem 4.2]{Carboni-Kelly-Pedicchio}. Every semi-abelian category is Mal'tsev \cite{Bourn1996, Borceux-Bourn}. Thus when $\Cc$ is semi-abelian, every simplicial object is Kan with respect to the class of regular projectives.

Note, however, that $\Cc$ need not have enough projectives for the internal Kan condition of~\cite{Carboni-Kelly-Pedicchio} to make sense. Following Remark~\ref{Remark-show-regular-epi}, the projective objects in the definition of the relative Kan property given here may be replaced by an enlargement of domain as in~\cite{Carboni-Kelly-Pedicchio}. In case there are enough projectives, of course the two notions do coincide.
\end{example}

\begin{example}\label{Group-Kan}
It is well known that the underlying simplicial set of a simplicial group is always Kan. This may be seen as a consequence of the previous example, because the category $\Gp$ is a Mal'tsev variety and the forgetful functor $U\colon{\Gp\to \Set}$ is represented by the group of integers $\Z$.
\end{example}

\subsection{Simplicial resolutions} Since $\Cc$ is an arbitrary category (without any extra structure) and $\Ac$ is just semi-abelian (rather than abelian), we have to be careful when considering simplicial resolutions of objects of $\Cc$. Definition~\ref{Definition-Resolution} seems to suit our purposes.

\begin{definition} 
Let $A=(A_n)_{n\geq-1}$ be an augmented simplicial object. A \defn{contraction} of $A$ is a family of maps $h_n\colon {A_n\to A_{n+1}}$, for $n\geq-1$, which satisfy $\del_0h_n=1_{A_n}$ and $\del_ih_n=h_{n-1}\del_{i-1}$ for $i>0$. A simplicial object that admits a contraction is called \defn{contractible}.
\end{definition}

\begin{definition}\label{Definition-Resolution}
Let $\Pc$ be a projective class. A \defn{$\Pc$-resolution of $X$} is an augmented simplicial object $A=(A_n)_{n\geq-1}$ with $A_{-1}=X$, where $A_n\in \Pc$ for $n\geq0$, and for every object $P\in \Pc$ the augmented simplicial set $\Hom(P,A)$ is Kan and contractible.
\end{definition}

In this paper, we focus on simplicial resolutions in a category $\Cc$ which are generated by a comonad $\G$ on $\Cc$. For any $\G$-projective object $P$, the simplicial set $\Hom(P,\G X)$ is contractible: choose a splitting $s$ for $\eps_{P}\colon{GP\to P}$; given a map $f\colon {P\to G^{n+1}X}$, define $h_n(f)=Gf \comp s\colon {P\to G^{n+2}X}$. The morphisms $h_n\colon {\Hom(P,G^{n+1}X)\to\Hom(P,G^{n+2}X)}$ then satisfy $\del_0  h_n=1_{\Hom(P,G^{n+1}X)}$ and $\del_i  h_n=h_{n-1}  \del_{i-1}$ for $i>0$. Thus they give a contraction of the simplicial set $\Hom(P,\G X)$. Later we assume that the category $\Cc$ and the projective class $\Pc$ generated by $\G$ are such that $\G X$ is Kan relative to $\Pc$ for any object $X$, so that $\G X$ is a $\Pc$-resolution of $X$.

In the case when $\Cc$ is a category with finite limits, there exists another definition of simplicial resolution, using \emph{simplicial kernels}. We give the definition of simplicial kernels here so that we can relate our Comparison Theorem of the next section with that of Tierney and Vogel \cite{Tierney-Vogel2}.

\begin{definition}\cite{Tierney-Vogel2}\label{Definition-Simplicial-Kernels}
Let 
$$(f_i\colon {X\to Y})_{0\leq i\leq n}$$
be a sequence of $n+1$ morphisms in the category $\Cc$. A \defn{simplicial kernel} of $(f_0,\ldots,f_n)$ is a sequence 
$$(k_i\colon {K\to X})_{0\leq i\leq n+1}$$
of $n+2$ morphisms in $\Cc$ satisfying $f_ik_j=f_{j-1}k_i$ for $0\leq i<j\leq n+1$, which is universal with respect to this property. In other words, it is the limit for a certain diagram in $\Cc$.
\end{definition}

For example, the simplicial kernel of one map is just its kernel pair.
If $\Cc$ has finite limits, simplicial kernels always exist. We can then factor any augmented simplicial object through its simplicial kernels as follows:
$$\xymatrix{{\cdots} \ar@<1.5ex>[rr] \ar@<.5ex>[rr] \ar@<-.5ex>[rr] \ar@<-1.5ex>[rr] \ar[dr] && A_2 \ar[rr] \ar@<1ex>[rr] \ar@<-1ex>[rr] \ar[dr] && A_1 \ar@<.5ex>[rr] \ar@<-.5ex>[rr] \ar[dr] && A_0 \ar[rr] && A_{-1}\\
& K_3 \ar@<1.5ex>[ur] \ar@<.5ex>[ur] \ar@<-.5ex>[ur] \ar@<-1.5ex>[ur] && K_2 \ar[ur] \ar@<1ex>[ur] \ar@<-1ex>[ur] && K_1 \ar@<.5ex>[ur] \ar@<-.5ex>[ur] & & }$$
Here the $K_{n+1}$ are the simplicial kernels of the maps $(\del_i)_{i}\colon {A_{n} \to A_{n-1}}$. This gives a definition of $\Pc$-exact simplicial objects:

\begin{definition}\cite{Tierney-Vogel2}
Let $\Pc$ be a projective class. An augmented simplicial object $A=(A_n)_{n\geq-1}$ is called \defn{$\Pc$-exact} when the comparison maps to the simplicial kernels and the map $\del_0\colon {A_0\to A_{-1}}$ are $\Pc$-epimorphisms. 
\end{definition}

\begin{remark}\label{Remark-exact-implies-contractible} It can be shown that for any $\Pc$-exact simplicial object $A$, the simplicial set $\Hom(P,A)$ is contractible for any $P\in \Pc$. We will call this property of $A$ \defn{relative contractability}.
\end{remark}

A resolution in the Tierney-Vogel sense is then a $\Pc$-exact augmented simplicial object $A$ in which all objects $A_n$ for $n\geq 0$ are in the projective class $\Pc$. For their definition they need the presence of simplicial kernels, so they have to assume for example that the category $\Cc$ has finite limits. In our definition all assumptions are on the comonad $\G$ or rather the induced projective class $\Pc$, and not on the category $\Cc$.  In the next section we will make clear the connections between our definition and theirs.

\section{The Comparison Theorem}\label{Section-Comparison-Theorem}
 
Let $\Pc$ be a projective class on $\Cc$.

\begin{lemma}\label{factormaps} 
Let $P\in\Pc$, and let $A$ be an augmented simplicial object for which the augmented simplicial set $\Hom(P,A)$ is contractible and Kan. Let $n\geq0$. Given a sequence of maps $(a_i\colon {P\to A_{n-1}})_{i\in [n]}$ satisfying $\del_i  a_j=\del_{j-1}  a_i$ for $i<j$, there is a map $a\colon {P\to A_n}$ with $\del_i  a=a_i$. 
\end{lemma}
\begin{proof}
Define the maps $b_{i+1}=h_{n-1}(a_i)$, where $(h_n)_{n\geq-1}$ is the contraction of the simplicial set $\Hom(P,A)$. These maps satisfy $\del_0  b_{i+1}=a_i$, and also $\del_j  b_{i+1}=\del_{i+1}  b_{j+1}$ for $i<j\leq n$, since $(\del_j  h_{n-1})(a_i)=h_{n-2}(\del_{j-1}  a_i)$, and $\del_{j-1}  a_i=\del_i  a_j$ for $i<j$. 
\begin{center}
\scalebox{.75}{\includegraphics{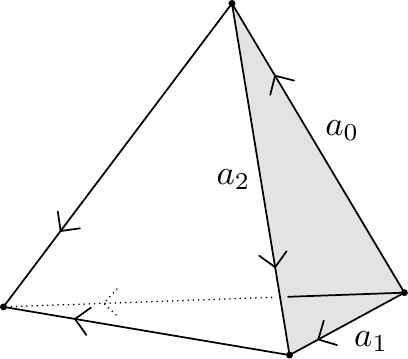}}
\end{center}
Thus they form an $(n+1,0)$-horn in the simplicial set $\Hom(P,A)$, and since we are assuming that this simplicial set is Kan, this horn has a filler $b\colon {P\to A_{n+1}}$. This gives the required map $a=\del_0b$. 
\end{proof}

\begin{remark}This lemma shows that in the presence of finite limits our $\Pc$-resolutions are also simplicial resolutions in the sense of Tierney and Vogel \cite{Tierney-Vogel2}; that is, the comparison maps to the simplicial kernels are $\Pc$-epimorphisms. Together with Remark~\ref{Remark-exact-implies-contractible} we see that $\Pc$-exactness and relative contractibility are equivalent in the situation when we have finite limits and any simplicial object is Kan relative to $\Pc$. So if $\Cc$ has finite limits, the Comparison Theorem 2.4 from \cite{Tierney-Vogel2} is more general than the one following in this section, but in the absence of finite limits our Comparison Theorem still works.\end{remark}

We now prove our Comparison Theorem using the above lemma.

\begin{theorem}[Comparison Theorem]\label{comparison}
Let $P$ be a simplicial object over $X$ with each $P_i\in\Pc$, and let $A$ be a simplicial object over $Y$, for which all the augmented simplicial sets $\Hom(P_i,A)$ are contractible and Kan. Then any map $f\colon {X\to Y}$ can be extended to a semi-simplicial map $f\colon {P\to A}$, and any two such extensions are simplicially homotopic.
\end{theorem}
\begin{proof}
We construct this semi-simplicial map inductively, using Lemma~\ref{factormaps}.
$$
\xymatrix{{\cdots} \ar@<1.5ex>[r] \ar@<.5ex>[r] \ar@<-.5ex>[r] \ar@<-1.5ex>[r] & P_{2} \ar@{.>}[d]^-{f_{2}} \ar[r] \ar@<1ex>[r] \ar@<-1ex>[r] & P_1 \ar@{.>}[d]^-{f_{1}} \ar@<.5ex>[r] \ar@<-.5ex>[r] & P_0 \ar@{.>}[d]^-{f_{0}} \ar[r] & X \ar[d]^-{f_{-1}=f} \\
{\cdots} \ar@<1.5ex>[r] \ar@<.5ex>[r] \ar@<-.5ex>[r] \ar@<-1.5ex>[r] & A_{2} \ar[r] \ar@<1ex>[r] \ar@<-1ex>[r] & A_1 \ar@<.5ex>[r] \ar@<-.5ex>[r] & A_0  \ar[r] & Y}
$$
Suppose the maps $f_j\colon {P_j\to A_j}$ are given for $-1\leq j< n$, and commute appropriately with the $\del_i$. This gives us $n+1$ maps $a_i\colon {P_n\to A_{n-1}}$, where $i\in[n]$, by composing the $\del_i\colon {P_n\to P_{n-1}}$ with $f_{n-1}$. These maps satisfy $\del_ia_j=\del_{j-1}a_i$ for $i<j$, since $\del_if_{n-1}=f_{n-2}\del_i$, and the $\del_j$ in $P$ satisfy the simplicial identities. Thus we can use Lemma~\ref{factormaps} to obtain the map $f_n\colon {P_n\to A_n}$ such that $\del_if_n=a_i=f_{n-1}\del_i$.

Now suppose $f\colon {P\to A}$ and $g\colon {P\to A}$ are two semi-simplicial maps commuting with $f\colon {X\to Y}$. We construct a  homotopy $h^n_i\colon {P_n\to A_{n+1}}$ for $n\geq 0$ and $0\leq i\leq n$, which satisfies $\del_0h^n_0=f_n$, $\del_{n+1}h^n_n=g_n$ and 
$$
\del_ih^n_j=\begin{cases}h^{n-1}_{j-1}\del_i& \text{ for } i<j\\
\del_ih^n_{i-1}& \text{ for } i=j\neq0\\
h^{n-1}_j\del_{i-1} & \text{ for } i>j+1.\end{cases}
$$
$h^0_0$ can be constructed using Lemma~\ref{factormaps}. Suppose the $h^k_j$ exist for $k<n$ and commute appropriately with the $\del_i$. Then $h^n_0$ must satisfy $\del_0h^n_0=f_n$, $\del_1h^n_0=\del_1h^n_1$ and $\del_ih^n_0=h^{n-1}_0\del_{i-1}$ for $i>1$. Of these maps, all are known except for $\del_1h^n_1$. Setting $a^0_0=f_n$ and $a^0_i=h^{n-1}_0\del_{i-1}$ for $i>1$, we form an $(n+1,1)$-horn in $\Hom(P_n,A)$. A filler for this horn gives $h^n_0$, and also $a^0_1=\del_1h^n_1$, which is needed for the next step. Now suppose $h^n_j$ are given for $j<l$, and we have $a^{l-1}_l=\del_lh^n_l=\del_lh^n_{l-1}$. Then $a^l_i=h^{n-1}_{l-1}\del_i$ for $i<l$, $a^l_l=a^{l-1}_l$ and $a^l_i=h^{n-1}_l\del_{i-1}$ for $i>l+1$ form an $(n+1,l+1)$-horn. A filler for this gives $h^n_l$ and $a^l_{l+1}=a^{l+1}_{l+1}$ for the next step. In the last step we have $a^n_i=h^{n-1}_{n-1}\del_i$ for $i<n$, $a^n_n=a^{n-1}_n=\del_nh^n_{n-1}$ and $a^n_{n+1}=g_n$. Then we use Lemma~\ref{factormaps} again to get $h^n_n$. 
\end{proof}

\section{Comonads generating the same projective class}\label{Section-Main-Theorem}
 
\subsection{Homotopy equivalence}
In this section we will assume that the category $\Cc$ and the projective class $\Pc$ generated by the comonad $\G$ are such that any augmented simplicial object $A$ which is relatively contractible is also Kan relative to $\Pc$. In particular the simplicial object $\G X$ is Kan relative to $\Pc$ for any object $X$. We will call such a projective class $\Pc$ a \defn{Kan projective class on $\Cc$}. Thus when $\G$ generates a Kan projective class, and $\K$ is a second comonad which generates the same projective class, the simplicial object $\K X$ is automatically also Kan relative to $\Pc$.

\begin{lemma}\label{GX-KX-homotopic}
Let $\G$ and $\K$ be two comonads on $\Cc$ which generate the same Kan projective class $\Pc$. Then the simplicial objects $\G X$~and $\K X$~are homotopically equivalent for any object~$X$.
\end{lemma}
\begin{proof}
Our assumptions on $\Cc$ and $\Pc$ imply that for any object $X$, the simplicial objects $\G X$ and $\K X$ are both $\Pc$-resolutions of $X$. Thus we can use the Comparison Theorem~\ref{comparison} to get semi-simplicial maps $f\colon {\G X\to \K X}$ and $g\colon {\K X\to \G X}$ which commute with the identity on $X$. 
\[
\xymatrix{{\cdots} \ar@<1.5ex>[r] \ar@<.5ex>[r] \ar@<-.5ex>[r] \ar@<-1.5ex>[r] & G^{3}X \ar[d]^-{f_{2}} \ar[r] \ar@<1ex>[r] \ar@<-1ex>[r] & G^{2}X \ar[d]^-{f_{1}} \ar@<.5ex>[r] \ar@<-.5ex>[r] & GX \ar[d]^-{f_{0}} \ar[r] & X \ar[d]^-{1_{X}} \\
{\cdots} \ar@<1.5ex>[r] \ar@<.5ex>[r] \ar@<-.5ex>[r] \ar@<-1.5ex>[r] & K^{3}X \ar[d]^-{g_{2}} \ar[r] \ar@<1ex>[r] \ar@<-1ex>[r] & K^{2}X \ar[d]^-{g_{1}}\ar@<.5ex>[r] \ar@<-.5ex>[r] & KX \ar[d]^-{g_{0}} \ar[r] & X \ar[d]^-{1_{X}} \\
{\cdots} \ar@<1.5ex>[r] \ar@<.5ex>[r] \ar@<-.5ex>[r] \ar@<-1.5ex>[r] & G^{3}X \ar[r] \ar@<1ex>[r] \ar@<-1ex>[r] & G^{2}X \ar@<.5ex>[r] \ar@<-.5ex>[r] & GX \ar[r] & X}
\]
Using the second part of the Comparison Theorem we see that both $fg$ and $gf$ are homotopic to the identity on $\K X$ and $\G X$ respectively. Thus $\G X$ and $\K X$ are homotopically equivalent. 
 \end{proof}
 
\begin{remark}
In this case we don't actually need the full strength of the second half of Lemma~\ref{comparison}. For any semi-simplicial map $f\colon {\G X\to \G X}$ which commutes with the identity on $X$, we can use the homotopy $h^n_i=(G^{i+1}f_{n-i})\sigma_i$ to see it is homotopic to the identity on $\G X$.
\end{remark}

\subsection{Simplicially homotopic maps in semi-abelian categories}
Given a functor $E\colon {\Cc\to \Ac}$, the simplicial objects $E\G X$ and $E\K X$ are still homotopically equivalent. We now show that, when $\Ac$ is a semi-abelian category, two simplicially homotopic semi-simplicial maps induce the same map on homology (see also \cite{VdLinden:Doc}). For this we need to define a special simplicial object, so that all the maps that form a simplicial homotopy are taken together to form a single semi-simplicial map. We do this by defining the following limit objects $A^I_n$.

\begin{notation}\label{Notation-A^J}

Suppose that $\Ac$ has finite limits and let $A$ be a simplicial object in $\Ac$. 
Put $A^{I}_{0}=A_{1}$ and, for $n> 0$, let $A^{I}_{n}$ be the limit (with projections $\pr_{1}$,\dots, $\pr_{n+1}\colon {A^{I}_{n}\to A_{n+1}}$) of the zigzag\index{zigzag}
$$\xymatrix@!@=.5cm{
A_{n+1}\ar[rd]_{\del_{1}} && A_{n+1} \ar[ld]^{\del_{1}} \ar[rd]_{\del_{2}}&& {} \ar@{}[d]|-{\displaystyle\cdots} \ar@{.>}[ld] \ar@{.>}[rd] && A_{n+1} \ar[ld]^{\del_{n}}\\
& A_{n} && A_{n} & {} & A_{n}}$$
in $\Ac$.

Let $\eps_{0} (A)_{n},\eps_{1} (A)_{n}\colon A^{I}_{n}\to A_{n}$ and $s (A)_{n}\colon A_{n}\to A^{I}_{n}$ denote the morphisms respectively defined by
$$\begin{aligned}
\eps_{0} (A)_{0} &= \del_{0}\\
\eps_{0} (A)_{n} &= \del_{0}  \pr_{1},
\end{aligned}
\quad 
\begin{aligned}
\eps_{1} (A)_{0} &= \del_{1}\\
\eps_{1} (A)_{n} &= \del_{n+1}  \pr_{n+1},
\end{aligned}
\quad\text{and}\quad 
s (A)_{n}= (\sigma_{0},\dots,\sigma_{n}).$$
\end{notation}

\begin{proposition}
Let $A$ be a simplicial object in a finitely complete category $\Ac$. Then the faces $\del_{i}^{I}\colon {A^{I}_{n}\to A^{I}_{n-1}}$ and degeneracies $\sigma_{i}^{I}\colon {A^{I}_{n}\to A^{I}_{n+1}}$ given by
\begin{gather*}
\begin{aligned}
\del_{0}^{I} &=\del_{0}  \pr_{2}\colon A^{I}_{1}\to A^{I}_{0}\\
\del_{1}^{I} &=\del_{2}  \pr_{1}\colon A^{I}_{1}\to A^{I}_{0}\\
\sigma_{0}^{I} &= (\sigma_{1},\sigma_{0})\colon A^{I}_{0}\to A^{I}_{1}
\end{aligned}
\qquad 
\pr_{j}  \del^{I}_{i} = \begin{cases}\del_{i+1}  \pr_{j} & \text{if $j\leq i$}\\
\del_{i}  \pr_{j+1} & \text{if $j>i$}  \end{cases}\colon A^{I}_{n}\to A_{n}\\
\pr_{k}  \sigma^{I}_{i} =
\begin{cases} \sigma_{i+1}  \pr_{k} & \text{if $k\leq i+1$}\\
\sigma_{i}  \pr_{k-1} & \text{if $k>i+1$} \end{cases}
\colon A^{I}_{n}\to A_{n+2},
\end{gather*}
for $i\in [n]$, $1\leq j\leq n$ and $1\leq k\leq n+2$, determine a simplicial object $A^{I}$. The morphisms mentioned in Notation~\ref{Notation-A^J} above form simplicial morphisms
$$\eps_{0} (A),\eps_{1} (A)\colon {A^{I}\to A}\qquad\text{and}\qquad s (A)\colon {A\to A^{I}}$$
such that $\eps_{0}(A)\comp s(A) = 1_{A} = \eps_{1}(A)\comp s(A)$. In other words, $(A^{I},\eps_{0}(A),\eps_{1}(A),s(A))$ forms a cocylinder on $A$.

Two semi-simplicial maps $f,g\colon{B\to A}$ are simplicially homotopic if and only they are homotopic with respect to the cocylinder $(A^{I},\eps_{0}(A),\eps_{1}(A),s(A))$: there exists a semi-simplicial map $h\colon{B\to A^{I}}$ satisfying $\eps_{0}(A)\comp h=f$ and $\eps_{1}(A)\comp h=g$.\noproof
\end{proposition}

Using a Kan property argument, we now give a direct proof that homotopic semi-simplicial maps have the same homology.

\begin{proposition}\label{Proposition-phi_0-Homology-Isomorphism}
Let $A$ be a simplicial object in a semi-abelian category $\Ac$; consider 
\[
\eps_{0} (A)\colon {A^{I}\to A}.
\]
For every $n\in \N$, $\H_{n}\eps_{0} (A)$ is an isomorphism.
\end{proposition}
\begin{proof}
Recall that in a semi-abelian category every simplicial object is Kan, relative to the class of regular epimorphisms. Using the Kan property, we show that the commutative diagram
$$
\xymatrix{\Nu_{n+1}A^{I} \ar@{-{ >>}}[rr]^-{\Nu_{n+1}\eps_{0} (A)} \ar[d]_-{d'_{n+1}} && \Nu_{n+1}A \ar[d]^-{d'_{n+1}}\\
\Zu_{n} A^{I} \ar@{-{ >>}}[rr]_-{\Zu_{n}\eps_{0} (A)} && {\Zu_{n}A}}
$$
is a regular pushout (see Definition~\ref{Definition-regepi-Pushout}); then it is a pushout, and Lemma~\ref{Lemma-pushout-implies-cokernels-iso} implies that the induced map $\H_{n}A^{I}\to \H_{n}A$ is an isomorphism.
Consider morphisms $z\colon Y_{0}\to \Zu_{n}A^{I}$ and $a\colon Y_{0}\to \Nu_{n+1}A$ that satisfy $d'_{n+1}\comp a=\Zu_{n}\eps _{0} (A)\comp z$.  It is enough to show that there exist a regular epimorphism $y\colon {Y\to Y_{0}}$ and a morphism $h\colon Y\to \Nu_{n+1}A^{I}$ satisfying $d'_{n+1}\comp h=z\comp y$ and $\Nu_{n+1}\eps_{0} (A)\comp h=a\comp y$: this implies that the comparison map to the pullback is a regular epimorphism, by Lemma~\ref{Lemma-show-regular-epi} and the fact that the morphisms of a limit cone form a jointly monic family.

We first sketch the geometric idea of this in the case $n=0$. Consider $a=a_{0}$ and $z=z_{0}$ as in Figure~\ref{Figure-Kan-Twice}; then (up to enlargement of domain) using the Kan property twice yields the needed $(h_{0},h_{1})$ in $\Nu_{1}A^{I}$.
\begin{figure}[htb]
\begin{center}
\scalebox{.75}{\includegraphics{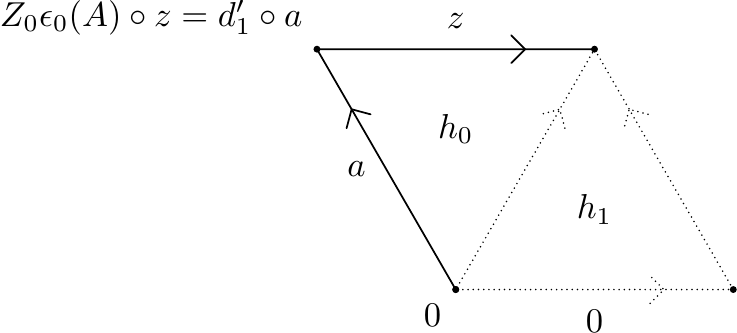}}
\end{center}
\caption{Using the Kan property twice to obtain $(h_{0},h_{1})$ in $\Nu_{1}A^{I}$.}\label{Figure-Kan-Twice}
\end{figure}

For arbitrary $n$, write
$$
a_{0}= \bigcap_{j} \ker \del_{j}\comp a\colon Y_{0}\to A_{n+1},
$$
and $(z_0,\dots,z_{n})= \bigcap_{j} \ker \del_{j}\comp z$. Note that as $z\colon Y_{0}\to \Zu_{n}A^{I}$, we have $\del^I_iz=0$ for $i\in[n]$, which implies $\del_iz_{j-1}=0$ for $i<j-1$ and $i>j$, where $1\leq j\leq n+1$. We also have $\del_jz_{j-1}=\del_jz_{j}$ for $1\leq j\leq n$ from the definition of the objects $A^I_n$. The map $a_0$ in turn satisfies $\del_ia_0=0$ for $i\in[n]$, and $\del_{n+1}a_0=\del_0z_0$. This last equality follows from $d'_{n+1}\comp a=\Zu_{n}\eps _{0} (A)\comp z$.

Suppose we have regular epimorphisms $y_k\colon {Y_k\to Y_{k-1}}$ for $1\leq k\leq n+2$, and morphisms $h_{k-1}\colon Y_{k}\to A_{n+2}$ satisfying $\del_{i}\comp h_{k-1}=0$ for $1\leq k\leq n+2$ and $i\notin \{k-1,k,n+2\}$, and $\del_{n+2}h_{k-1}=z_{k-1}y_1\cdots y_{k}$ for $1\leq k\leq n+1$, and also $\del_0h_0=a_0y_1$. We set $y=y_1\comp \cdots \comp y_{n+2}$.
This gives us the required map
$$
h=(h_{0}\comp y_{2}\comp \cdots \comp y_{n+2},\dots, h_{n+1})\colon Y_{n+2}\to \Nu_{n+1}A^{I}
$$
which satisfies $d'_{n+1}\comp h=z\comp y$ and $\Nu_{n+1}\eps_{0} (A)\comp h=a\comp y$.

We construct these maps $y_k$ and $h_k$ inductively. To get $h_0$ we form an $(n+2,1)$-horn $(b_i\colon {Y_0\to A_{n+1}})_{i}$ of $A$ by setting $b_0=a_0$, $b_{n+2}=z_0$ and $b_i=0$ for $1<i<n+2$. A filler for this horn gives $y_1\colon {Y_1\to Y_0}$ and $h_0\colon {Y_1\to A_{n+2}}$ satisfying $\del_0h_0=a_0y_1$. Now suppose for $1\leq k\leq n+1$ we have $a_{k}=\del_kh_{k-1}\colon Y_{k}\to A_{n+1}$ with $\del_{i}\comp a_{k}=0$ for $i\in [n]$, and $\del_{n+1}a_k=\del_kz_{k-1}y_1y_2\cdots y_k$. Then we can form an $(n+2,k+1)$-horn by setting $b_{k}=a_{k}$, $b_{n+2}=z_{k}y_1\cdots y_k$ and $b_i=0$ for $i<k$ and $k+1<i<n+2$, which induces $y_{k+1}$ and $h_{k}$ with the desired properties. 
\end{proof}

\begin{remark}
A homology functor $\H_{n}$ involves an implicit choice of colimits: the cokernels involved in the construction of the $\H_{n}A$. We may, and from now on we will, assume that these colimits are chosen in such a way that $\H_{n}\eps_{0}(A)$ is an \emph{identity} instead of just an isomorphism. This gives us the equality in the next corollary.
\end{remark}

\begin{corollary}\label{Corollary-Homotopic-then-Same-Homology}
If $f\simeq g$ then, for any $n\in \N$, $\H_{n}f= \H_{n}g$.
\end{corollary}
\begin{proof}
Proposition~\ref{Proposition-phi_0-Homology-Isomorphism} states that $\H_{n}\eps_{0}(A)$ is an isomorphism; by a careful choice of colimits in the definition of $\H_{n}$, we may assume that $\H_{n}\eps_{0}(A)=1_{\H_{n}A}=\H_{n}\eps_{1}(A)$.
\[
\xymatrix{&&&&A \ar@{=}[d]\\
B \ar[rr]^-h \ar@/^/[rrrru]^-g \ar@/_/[rrrrd]_-f &&A^I \ar[rru]_-{\eps_1(A)} \ar[rrd]^-{\eps_0(A)}&&A\ar[ll]|{s(A)}\\
&&&&A\ar@{=}[u]}
\]
If now $h$ is a homotopy $f\simeq g$, then $\H_{n}f=\H_{n}\epsilon_{0}(A)\comp \H_{n}h=\H_{n}\epsilon_{1}(A)\comp \H_{n}h=\H_{n}g$.
\end{proof}

\subsection{Isomorphism between the homology functors}
Using the above, we can now prove our Main Theorem.

\begin{theorem}\label{Main-Theorem}
Let $\G$ and $\K$ be two comonads on $\Cc$ which generate the same Kan projective class $\Pc$. Let $E\colon {\Cc\to\Ac}$ be a functor into a semi-abelian category. Then the functors $\H_n(-,E)_\G$ and $\H_n(-,E)_\K$ from $\Cc$ to $\Ac$ are isomorphic for all $n\geq1$.
\end{theorem}
\begin{proof}
It follows from Lemma~\ref{GX-KX-homotopic} that the simplicial objects $E\G X$ and $E\K X$ are homotopically equivalent. Thus Corollary~\ref{Corollary-Homotopic-then-Same-Homology} implies that $\H_n(X,E)_\G\iso \H_n(X,E)_\K$. Given a map $f\colon {X\to Y}$, the two semi-simplicial maps
\[
\xymatrix{\G X\ar[r] &\K X \ar[r]^-{\K f} & \K Y} \qquad\text{and}\qquad\xymatrix{\G X\ar[r]^-{\G f} & \G Y \ar[r] & \K Y}
\]
are both semi-simplicial extensions of $f$, so they are homotopic by the Comparison Theorem~\ref{comparison}. Again using Corollary~\ref{Corollary-Homotopic-then-Same-Homology}, we see that the square
\[
\xymatrix@=1cm{\H_{n}\G X \ar[r]^-{\H_{n}\G f} \ar[d]_-{\cong} & \H_{n}\G Y \ar[d]^-{\cong} \\
\H_{n}\K X \ar[r]_-{\H_{n}\K f} & \H_{n}\K Y}
\]
commutes, which proves that the isomorphisms are natural.
\end{proof}

\begin{remark}In fact, the above isomorphism is also natural in the second variable. If $\alpha\colon {E\To F}$ is a natural transformation, then the square
$$\xymatrix@=1cm{\H_{n}(X,E)_\G \ar[r]^-{\iso} \ar[d]_-{\H_{n}(X,\alpha)_\G} & \H_{n}(X,E)_\K \ar[d]^-{\H_{n}(X,\alpha)_\G} \\
\H_{n}(X,F)_\G \ar[r]_-{\iso} & \H_{n}(X,F)_\K}$$
also commutes, since 
$$\xymatrix@=1cm{E\G X\ar[r] \ar[d]_-{\alpha_{\G X}} & E\K X \ar[d]^-{\alpha_{\K X}} \\
F\G X \ar[r] & F\K X}$$
already commutes. 
\end{remark}

\begin{remark} We could define homology just using a projective class instead of a comonad, since the Comparison Theorem and Corollary~\ref{Corollary-Homotopic-then-Same-Homology} imply that \emph{any} $\Pc$-resolution of $X$ will give the same homology. Consider for example the following (Tierney-Vogel) resolution in a category with finite limits:

Given an object $X$, there is a $\Pc$-projective object $X_0$ with a $\Pc$-epimorphism $\del_0\colon {X_0 \to X}$, since there are enough $\Pc$-projectives. We can call this a presentation of $X$. Take the kernel pair of $\del_0$, and take the presentation of the resulting object to get $X_1$. Composition gives two maps $\del_0$ and $\del_1$ to $X_0$, and we can take the simplicial kernel of these and the presentation of the resulting object to get $X_3$ and so on. This gives a resolution in the Tierney-Vogel sense \cite{Tierney-Vogel2}. When $\Pc$ is a Kan projective class, it is also a $\Pc$-resolution in our sense and thus gives the same homology. This resolution is often easier to work with than the functorial $\G X$.
\end{remark}

\section{Examples}\label{Section-Examples}

We first give some examples of categories $\Cc$ and comonads $\G$ which generate a Kan projective class, and then give a specific example in the category of $R$-modules of two comonads giving the same projective class. 

\subsection{Some valid contexts}
If $\Cc$ is an additive category, our condition on the projective class $\Pc$ is easily satisfied by any comonad on $\Cc$, since then for any simplicial object $A$ and any object $P$, the simplicial set $\Hom(P,A)$ is actually a simplicial group and thus Kan (cf.\ Example~\ref{Group-Kan}). This includes many of the examples in Barr and Beck's paper \cite{Barr-Beck}, for example the comonad $\G$ on the category $R\text{-}\Mod$ of (left) $R$-modules generated by the forgetful/free adjunction to $\Set$, the forgetful/free comonad on the category $\Comm$ of commutative rings, and also the comonad on the category $K\text{-}\Alg$ of associative $K$-algebras generated by the forgetful functor to $K\text{-}\Mod$. See \cite{Barr-Beck} for more additive examples. 

When $\Cc$ is a regular category and the projective class $\Pc$ is the class of regular projectives, as remarked in Section~\ref{Section-Relative-Kan-Property} saying that a simplicial object $A$ is Kan relative to $\Pc$ is the same as saying $A$ is internally Kan in $\Cc$. Thus when $\Cc$ is also Mal'tsev, every simplicial object is Kan \cite[Theorem 4.2]{Carboni-Kelly-Pedicchio}, and $\G X$ is a $\Pc$-resolution. This includes the forgetful/free comonads on the category $\Gp$ of groups, $\Rng$ of non-unital rings, $\XMod$ of crossed modules, etc. 

Given a comonad $\G$ on a category $\Cc$ which comes from an adjunction 
$$\xymatrix{\Cc\ar[rd]_{U} \ar[rr]^{G}&& \Cc\\
&\Dc \ar[ur]_{F}&}$$
we can determine the class of morphisms of the projective class $(\Pc,\Ec)$ generated by $\G$ in the following way:

Given an object $A$ and a morphism $e\colon {B\to C}$ in $\Cc$, the diagram 
$$\xymatrix{&GA\ar[d]^{f}\\ 
B\ar[r]^e&C}$$
corresponds via the adjunction to 
$$\xymatrix{&UA\ar[d]\\
UB\ar[r]^{Ue}&UC}$$
If $e$ is in $\Ec$, by choosing $A=C$ and $f=\eps_C$, we see that $Ue$ must be split in $\Dc$, since $\eps_C$ corresponds under the adjunction to $1_{UC}$. Conversely if $Ue$ is a split epimorphism in $\Dc$, we can factor any map $UA\to UC$ over $Ue$, which implies that we can factor any map $f\colon {GA\to C}$ over $e$ in $\Cc$, thus $e\in \Ec$. Therefore the class $\Ec$ is exactly the class of morphisms whose images under $U$ are split in $\Dc$. Thus when $\Cc$ is a variety and $U$ is the forgetful functor to $\Set$, we will always get the class of regular projectives on $\Cc$.

\subsection{Two comonads on $R\text{-}\Mod$}
Given a ring homomorphism $\phi\colon {S\to R}$, every $R$-module can also be considered as an $S$-module by restricting the $R$-action to $S$ via $ \phi$. This gives rise to an adjunction $R\otimes_S(-)\dashv  \Hom_R(R,-)$ between the categories of modules, where $R$ is viewed as an $S$-module.
 
We now consider the following situation: 
Let $S_1$ and $S_2$ be Morita-equivalent rings, where the equivalence between $S_1\text{-}\Mod$ and $S_2\text{-}\Mod$ is induced by a ring homomorphism $\psi\colon {S_1\to S_2}$. Let $R$ be another ring, with ring homomorphisms as below which make the diagram commute:
$$\xymatrix@=.5cm{S_1\ar[rrd]^{\phi_1} \ar[dd]^{\psi}&&\\
&&R\\
S_2 \ar[rru]^{\phi_2}&&}$$

For each $i=1,2$ this gives us a comonad on $R\text{-}\Mod$ using the adjunction above:
$$\xymatrix{R\text{-}\Mod \ar[rd]_{U_i} \ar[rr]^{G_i}&& R\text{-}\Mod \\
&S_i\text{-}\Mod \ar[ur]_{R\otimes_{S_i}(-)}&}$$

We write $U_i$ for the forgetful functor $\Hom_R(R,-)$ from $R$-modules to $S_i$-modules.

As seen above the projective class generated by $\G_i$ is given by the class of maps in $R\text{-}\Mod$ which are split as $S_i$-module maps. Since $\psi$ induces an equivalence between $S_1\text{-}\Mod$ and $S_2\text{-}\Mod$, any $R$-module map $e$ has the property that $U_1(e)$ is split if and only if $U_2(e)$ is split. Thus the comonads $\G_1$ and $\G_2$ induce the same projective class, and thus give rise to the same homology on $R\text{-}\Mod$. As mentioned in \cite{Barr-Beck} this homology is Hochschild's $S$-relative $\Tor$; so we see we get the same relative $\Tor$ functor for two Morita equivalent rings when the equivalence is induced by a ring homomorphism.

\section{Acknowledgements}
We would like to thank Tomas Everaert, Alexander Frolkin, Martin Hyland, Peter Johnstone and Alexander Shannon for useful comments and suggestions.

\end{document}